\newif\ifarxiv
\arxivfalse

\documentclass[twocolumn]{autart}    
\usepackage{algorithmicx}
\usepackage{algorithm}
\usepackage[noend]{algpseudocode}
\usepackage{subcaption}

\usepackage{amsmath}
\usepackage{amssymb}
\usepackage{import}
\usepackage{color}
\usepackage{ulem}
\usepackage{multirow}
\usepackage{cite} 
\usepackage{graphicx}          

\usepackage{array}
\usepackage{comment}
\usepackage{todonotes}
\usepackage{xspace} 

\newcommand{\inlinetodo}[1]{\todo[inline,size=\scriptsize]{#1}}
\newcommand{\argmin}{\operatornamewithlimits{argmin}}
\newcommand{\argmax}{\operatornamewithlimits{argmax}}
\newcommand{\Studentst}{Student's~$t$\xspace}
\newcommand{\hip}{\text{hip}}

\newcommand{\mh}{m_h}
\newcommand{\GD}{\text{G}\!\downarrow}
\newcommand{\GR}{\text{G}\!\uparrow}
\newcommand{\AD}{\text{A}\!\downarrow}
\newcommand{\AR}{\text{A}\!\uparrow}
\newcommand{\A}{\text{A}}
\newcommand{\G}{\text{G}}
\newcommand{\state}{state\xspace}
\newcommand{\states}{states\xspace}
\newcommand{\Proc}{\mathcal{F}}
\newcommand{\Meas}{\mathcal{H}}

\newcommand{\pnoise}{\sigma}
\newcommand{\mnoise}{\delta}
\newcommand{\cstatet}{x_t}
\newcommand{\cstate}{x}
\newcommand{\cstatetp}{x_{t+1}}
\newcommand{\obst}{y_t}

\newcommand{\delem}{w}
\newcommand{\ddim}{M}
\newcommand{\cdim}{n}
\newcommand{\odim}{d}
\newcommand{\Indicator}{I}
\newcommand{\wtrue}{w}
\newcommand{\rwest}{\widetilde{w}}
\newcommand{\west}{\widehat{w}}

\newcommand{\set}[1]{\left\{ #1 \right\}}

\newcommand{\R}{\mathbb{R}}

\newcommand{\into}{\rightarrow}

\newtheorem{theorem}{Theorem}
\newtheorem{example}{Example}
\newtheorem{lemma}{Lemma}
\newtheorem{corollary}{Corollary}

\newtheorem{proposition}{Proposition}
\newtheorem{assumption}{Assumption}
\newtheorem{definition}{Definition}

\newtheorem{remark}{Remark}
\def\BT{\begin{theorem}}
\def\ET{\end{theorem}}
\def\BL{\begin{lemma}}
\def\EL{\end{lemma}}
\def\BP{\begin{proposition}}
\def\EP{\end{proposition}}
\def\BC{\begin{corollary}}
\def\EC{\end{corollary}}
\def\BD{\begin{definition}}
\def\ED{\end{definition}}
\def\BA{\begin{assumption}}
\def\EA{\end{assumption}}
\def\BR{\begin{remark}}
\def\ER{\end{remark}}
\def\BE{\begin{example}}
\def\EE{\end{example}}

\newcommand*\Let[2]{\State #1 $\gets$ #2}
\newcommand{\figref}[1]{Figure~\ref{#1}}

\begin{document}

\begin{frontmatter}

\title{Offline state estimation for hybrid systems \\ via nonsmooth variable projection} 

\thanks[footnoteinfo]{%
This material is based on work supported in part by the Washington Research Foundation Data Science Professorship, in part by the U. S. Army Research Laboratory and the U. S. Army Research Office under grant \#W911NF-16-1-0158, and in part by the National Science Foundation Cyber-Physical Systems program award \#1565529.
}

\author[AM]{Jize Zhang}\ead{jizez@uw.edu},    
\author[ECE]{Andrew M. Pace}\ead{apace2@uw.edu},               
\author[ECE]{Samuel A. Burden}\ead{sburden@uw.edu},  
\author[AM]{Aleksandr Aravkin}\ead{saravkin@uw.edu}

\address[AM]{Department of Applied Mathematics, University of Washington, Seattle}  
\address[ECE]{Department of Electrical \& Computer Engineering, University of Washington, Seattle}             


\begin{abstract}                          
%
{A \emph{hybrid} dynamical system switches between dynamic regimes at time- or state-triggered events.}
We propose an offline algorithm that simultaneously estimates discrete and continuous components of a hybrid system's state.
We formulate state estimation as a continuous optimization problem by relaxing the discrete component 
and use a robust loss function to accommodate large changes in the continuous component during switching events.
Subsequently, we develop a novel nonsmooth variable projection algorithm with Gauss-Newton updates to solve the state estimation problem and prove the algorithm's global convergence to stationary points.
We demonstrate the effectiveness of our approach on simple piecewise-linear and -nonlinear mechanical systems undergoing intermittent impact.

\end{abstract}

\end{frontmatter}

\section{Introduction}

%
This paper considers the problem of using noisy measurements from a piecewise-continuous trajectory 
to estimate a hybrid system's state.
%
The state estimation problem has been extensively studied in \emph{classical} dynamical systems whose states evolve 
according to one 
(possibly time--varying) 
smooth model.
This problem is fundamentally more challenging for hybrid systems 
since the set of discrete state%
\footnote{%
We refer to the discrete component of the hybrid system state as the \emph{discrete state}, and similarly refer to the \emph{continuous state}, although the state of the hybrid system is specified by both the discrete and continuous components.
}
sequences generally grows combinatorially in time.

When the discrete state sequence and switching times are known \emph{a priori} or directly measured, only the continuous state needs to be estimated,
yielding a classical state estimation problem;
this approach has been applied to
piecewise-linear systems~\cite[Chap.~4.5]{Stengel1994OptimalControlEstimation}
and to
nonlinear mechanical systems undergoing impacts~\cite{MeniniTornambe2001AsymptoticTrackingPeriodica}.
When the discrete state
is not known or measured, 
estimating both the discrete and continuous states simultaneously improves estimation performance.
%
One approach uses a bank of filters, each tuned to one discrete state, 
and selects the discrete states as the filter with the lowest residual \cite[\S4.1]{BalluchiEtAl2002DesignObserversHybrid}.
This filter bank method 
has been applied to 
hybrid systems with 
linear dynamics \cite[\S4.1]{BalluchiEtAl2003ObservabilityHybridSystems} \cite{Gomez-GutierrezEtAl2011SlidingModeObserver}, 
nonlinear dynamics~\cite{BarhoumiEtAl2012ObserverDesignClassesa},
%
and jumps in the continuous state when the discrete state changes \cite{BalluchiEtAl2013DesignDynamicalObservers}.
%
%
Likewise, particle filter methods for hybrid systems \cite{BlomBloem2007ExactBayesianParticle,DoucetEtAl2001ParticleFiltersState,Seah2009-ev}
use a collection of filters, identified as particles, and are applicable to more general nonlinear process dynamics.
Particle filters and filter banks are effective when the number of discrete states and dimension of continuous state spaces are small.

Another approach formulates a moving-horizon estimator
over both the continuous and discrete \states,
resulting in a mixed-integer optimization problem~\cite{BemporadEtAl1999MovingHorizonEstimation}.
The inherently discrete nature of the problem formulation
enables estimation of the exact sample when the discrete state switches, 
at the expensive of combinatorial growth of the set of discrete decision variables as the horizon increases.
Multiple methods have been developed to mitigate the challenge posed by this combinatorial complexity. 
One approach entails 
summarizing past measurements and state estimates with a penalty term in the the objective function~\cite{Ferrari-TrecateEtAl2002MovingHorizonEstimation}.
Another approach, applicable to systems with bounded noise, entails 
restricting the set of possible discrete state sequences 
using \emph{a priori} knowledge of the system~\cite{AlessandriEtAl2005RecedinghorizonEstimationSwitching,AlessandriEtAl2007MinimumDistanceRecedingHorizonState}.

An alternative approach to circumventing the combinatorial challenge entailed by exactly estimating the discrete state sequence involves \emph{relaxing} the discrete state estimate to take on continuous values
as in~\cite{BakoLecoeuche2013SparseOptimizationApproach,Johnson2016ObservabilityObserverDesign}.
The latter reference uses a sparsity-promoting convex program whose objective incorporates a nonsmooth penalty across all possible discrete state sequences, and guarantees the estimate converges to the true continuous and discrete states.
Both approaches are formulated for piecewise-linear systems whose continuous states do not jump when switching between subsystems; in the language of hybrid systems, the continuous states are \emph{reset} using the identity function.
%

\subsection*{Our approach and contributions}

We propose an offline algorithm for estimating the state of hybrid systems with nonlinear dynamics, non--identity resets, and noisy process and observation models. 
Our starting point is the optimization perspective on generalized and robust state estimation~\cite{AravkinEtAl2016GeneralizedKalmanSmoothing,AravkinEtAl2012RobustTrendfollowingKalman}. 
To formulate state estimation as a continuous optimization problem, 
we relax the discrete state to take on continuous values as in prior work. 
Unlike prior work, we model process noise using the \Studentst distribution,
which allows large innovations and makes the method applicable to systems with non-identity resets.

In combination, these elements yield a nonsmooth nonconvex continuous optimization formulation (Sec.~\ref{sec:problem}). 
We develop a Gauss-Newton type algorithm to solve this problem 
%
and prove the algorithm globally converges to stationary points (Sec.~\ref{sec:algorithm}). 
The problem formulation and algorithm is evaluated on piecewise-linear and -nonlinear hybrid system models (Sec.~\ref{sec:experiments}).

\section{Problem formulation}
\label{sec:problem}
\subsection{Process and observation models}

We consider a class of \emph{discrete-time switched systems}
\begin{equation}
\begin{aligned}
x_{t+1} &= \sum_{m = 1}^{\ddim} \Proc_m(x_t) w_t[m] + \pnoise_t\\
y_t & = \Meas_t( x_t )+ \mnoise_t
 \label{eq:hybrid_dynamics}
 \end{aligned}
 \end{equation}
%
where
$m\in\ddim$ indexes the continuously-differentiable process model
$\Proc_m\colon\mathbb{R}^\cdim\into\mathbb{R}^\cdim$,
$\Meas_t\colon  \mathbb{R}^\cdim \into \mathbb{R}^\odim$ is the continuously-differentiable measurement model
that generates 
observations $\obst\in\mathbb{R}^\odim$ 
of the hidden continuous state $\cstatet\in\mathbb{R}^\cdim$,
$\pnoise_t, \mnoise_t$ are process and measurement noises,
and
$w_t \in \mathcal{D}^\ddim$ is a \emph{one-hot} vector%
\footnote{
$w\in \R^\ddim$ is \emph{one-hot} if 
$w[i] \in \{0,1\}$
for all $i\in\set{1,\dots,\ddim}$
and
$1^Tw = 1$;
$\mathcal{D}^\ddim\subset\R^\ddim$
denotes the set of one-hot vectors.}%
that indicates which process model is
active at time $t$.
Note that the observations do not depend explicitly on the active model $m$, which must be inferred from measurements of the continuous state $\cstatet$.

The model $m$ that is active during each time step 
may be 
determined by an exogenous signal, prescribed as a function of time or state, or some combination thereof.
Thus, the equation in~\eqref{eq:hybrid_dynamics} can represent the process and observation models of a wide variety of hybrid systems.
We are motivated theoretically and experimentally to focus on cases where the active model is constant for many time steps, only occasionally switching to a new model.
Such cases arise, for instance,
when \emph{sampling} trajectories of a hybrid dynamical system with a fixed timestep
(so long as the timestep is much smaller then the average \emph{dwell time}~\cite{HespanhaMorse1999StabilitySwitchedSystemsa}).

When the discrete state changes in a hybrid system, the continuous state may change abruptly according to a \emph{reset} map.
As an example, the velocity of a rigid mass changes abruptly when it impacts a rigid surface~\cite{Lotstedt1982MechanicalSystemsRigid}.
%
%
Empirically, these discrete reset dynamics are much more poorly characterized than their continuous counterparts. 
For instance, whereas the ballistic trajectory of a rigid mass is well-approximated by Newton's laws,
the abrupt change in velocity that occurs at impact is not consistent with any established impact law~\cite{FazeliEtAl2017LearningDataEfficientRigidBody}.
Including such a reset in the system model~\eqref{eq:hybrid_dynamics} will introduce bias 
into the state estimate
because the model will generate erroneous predictions at resets, diminishing the accuracy of estimated states at nearby times.
This observation motivates us in the next section to account for the effect of unknown resets as part of the process noise.

\subsection{Process noise and observation noise models}
\label{sec:noise_models}

Instead of incorporating continuous state resets explicitly into the model~\eqref{eq:hybrid_dynamics}, we introduce a distributional assumption on the process noise $\pnoise_t$ that accepts large 
instantaneous
changes in the continuous state estimate.  
Specifically, we assume that process noise $\pnoise_t$ follows a \Studentst distribution.
Compared with the commonly-used Gaussian distribution, the \emph{heavy-tailed} \Studentst is tolerant to large deviations in the estimate of the hidden continuous state $x_t$~\cite{aravkin2014robust}. 
%
Hence, the \Studentst error model 
allows an instantaneous change in the state that is consistent with~\eqref{eq:hybrid_dynamics} before and after the change.
The negative log-likelihood of the \Studentst (as a function of $\pnoise_t$) is given by 
\[
r\log\left( r + \left\|Q^{-1/2}\pnoise_t\right\|^2\right) - r\log(r),
\]
where $r$ is the degrees of freedom parameter of the \Studentst, and $Q$ is the covariance matrix.

If the continuous state $\cstate_t$ was known, then any residual between the predicted observations $\Meas_t(\cstate_t)$ and 
actual measurements $\obst$ at time $t$ is due to measurement noise;
in particular, the residual does not exhibit large deviations due to continuous state resets at switching times.
Thus, we assume the measurement noise $\mnoise_t$ follows the usual Gaussian distribution, with negative log-likelihood 
\[
\frac{1}{2}\left\| R^{-1/2} \mnoise_t\right\|^2,
\]
where $R$ is the covariance matrix.
The plots below provide a comparison between the  
probability density (left) and the negative log-likelihood (right)
for the
scalar
Gaussian (solid blue) 
and 
\Studentst distributions (dashed red; degree-of-freedom $r=1$).
\begin{tabular}{c c}
\includegraphics[height=.75in]{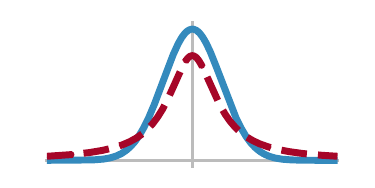} &
\includegraphics[height=.75in]{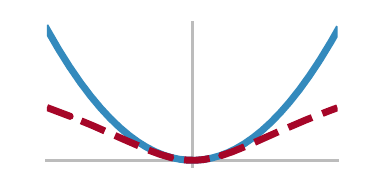}\\
probability density & negative log-likelihood
\end{tabular}

\subsection{State estimation problem formulation}
The \emph{maximum a posteriori} (MAP) likelihood for estimating states of~\eqref{eq:hybrid_dynamics}  
is given by 
\begin{equation}
\label{eq:MAP}
\begin{aligned}
&\min_{x_t,w_t\in\mathcal{D}^m} \sum_{t=0}^{T-1}
 \frac{1}{2}\left\| R^{-1/2}\left(\obst - \Meas_t(\cstatet) \right)\right\|^2 + \\
& r\log\left(
r + \left\|Q^{-1/2}\left(\cstatetp -\sum_{m = 1}^{M}\Proc_m(\cstatet)w_{t}[m]\right) \right\|^2\right).
\end{aligned}
\end{equation}
Problem~\eqref{eq:MAP} is a nonlinear mixed-integer program with respect to both the
continuous ($\cstatet$) 
and 
discrete ($w_t$)
decision variables. 
We can significantly simplify the structure by establishing the following lemma. 



\begin{lemma}[Formulation Equivalence]
Given $w \in \mathcal{D}^\ddim$, any vectors $x_1, x_2$, models $\Proc_i$, 
and any penalty functional $g$,  we have 
\begin{align*}
& \min_{w\in \mathcal{D}^\ddim} \quad  g\left (x_2 - \sum_{m=1}^\ddim w[m]\Proc_m(x_1)\right) \\
= &\min_{w\in\mathcal{D}^\ddim} \quad\sum_{m=1}^\ddim w[m] g\left(x_2- \Proc_m(x_1)\right)\\
\intertext{and}
& \argmin_{w\in \mathcal{D}^\ddim}  g\left(x_2 - \sum_{m=1}^\ddim w[m]\Proc_m(x_1)\right) \\
 = & 
\argmin_{w\in\mathcal{D}^\ddim} \quad \sum_{m=1}^\ddim w[m] g\left(x_2 - \Proc_m(x_1)\right).
\end{align*}
\end{lemma}
{\bf Proof}:

Since $w \in \mathcal{D}^\ddim$ for both problems, there are only $m$ possible values for both objective functions, i.e.
\[ g(x_2-\Proc_1(x_1)), \quad g(x_2-\Proc_2(x_1)),\quad \dots, \quad g(x_2-\Proc_m(x_1)).\]
Hence, the minimum objective value for both problems will be $\min_i g(x_2-\Proc_i(x_1))$ and every minimizer is a 
one-hot vector that selects a minimum value. $\square$

Therefore, an equivalent formulation to~\eqref{eq:MAP} is given by 
\begin{equation}
\label{eq:MAPalt}
\begin{aligned}
& \min_{x_t,w_t\in\mathcal{D}^\ddim} \sum_{t=0}^{T-1} \bigg (
 \frac{1}{2}\left\| R^{-1/2}\left(\obst - \Meas_t(\cstatet) \right)\right\|^2 + \\
& \sum_{m=1}^\ddim w_{t}[m]r\log\left( r + \left\|Q^{-1/2}\left(x_{t+1} -\Proc_m(\cstatet)\right) \right\|^2\right)\bigg).
\end{aligned}
\end{equation}
Although still a mixed-integer program, this reformulation exhibits linear coupling between the discrete variables $w_t$ and continuous variables $\cstatet$. 
We will leverage this linear coupling when we develop our estimation algorithm based on the relaxed problem formulation introduced in the next section.
  
\subsection{Relaxed state estimation problem formulation}
\label{sec:relax_smooth}
Ultimately, the discrete state estimate will be specified as a one-hot vector,
$w_t \in \mathcal{D}^\ddim \subset \R^\ddim$.
To formulate a continuous optimization problem that approximates the mixed-integer problem formulated in the previous section,  
we relax the decision variable $w_t$ to take values in the 
convex hull $\Delta^\ddim$ of $\mathcal{D}^\ddim$.%
\footnote{
We use
${\Delta^\ddim:= \{ w \in [0,1]^\ddim: 1^Tw = 1\}}$ to denote the simplex in $\mathbb{R}^\ddim$. 
}
%
%
The optimal relaxed $w_t$ will generally lie on the interior of the simplex, so we project the result from our relaxed optimization problem 
to return the one-hot discrete state estimate.
Since this relaxation-optimization-projection process tends to induce frequent changes in the discrete state estimate, 
we introduce a smoothing term on $w_t$,
\[
\nu\|w_{t+1} - w_t\|_2^2,
\]
yielding the continuous relaxation of~\eqref{eq:MAPalt} given by 
 \begin{equation}
 \label{eq:ss_full_obj}
\begin{aligned}
 & \min_{x_t,w_t}  f(x,w):=  \sum_{t=0}^{T-1}\bigg(
 \frac{1}{2}\left\| R^{-1/2}\left(\obst - \Meas_t(\cstatet) \right)\right\|^2 \\
 & +  \sum_{m=1}^\ddim w_{t}[m]r\log\left( r + \left\|Q^{-1/2}\left(\cstatetp -\Proc_m(\cstatet)\right) \right\|^2\right) \\
& + \nu\|w_{t+1} - w_t\|_2^2 + \Indicator(w_t | \Delta^\ddim)\bigg),
\end{aligned}
 \end{equation}
 where $x$ is the concatenated variable containing all $x_t$, 
 $w$ is the concatenated variable containing all $w_t$, and 
 \emph{convex indicator function} $\Indicator$ is defined by 
 \[
 \Indicator(w_t | \Delta^\ddim) := 
 \begin{cases} 0, &  w_t \in \Delta^\ddim; \\ \infty, & w_t \not\in\Delta^\ddim. 
 \end{cases}
 \]
The optimal relaxed 
discrete state estimate $w_t \in \Delta^\ddim$
is projected onto $\mathcal{D}^\ddim$ by choosing the (unique) one-hot vector whose $\argmax_i {w}_t[i]$ component is equal to~1.

\section{State estimation algorithm}
\label{sec:algorithm}

We develop an algorithm to solve the relaxed state estimation problem formulated in~\eqref{eq:ss_full_obj} using two key ideas:

\begin{enumerate}
\item nonsmooth variable projection;
\item Gauss-Newton descent with \Studentst penalties.
\end{enumerate}

These two ideas are explained in the next two subsections, followed by a convergence analysis in the third subsection.

\subsection{Nonsmooth variable projection}

The first idea is to pass to the \emph{value function}, projecting out (partially minimizing over) the $w$ variables. 
Define 
\begin{equation}
\label{eq:vf}
\begin{aligned}
v(x) := \min_{w} f(x,w)
\end{aligned}
\end{equation}
with $f(x,w)$ as in~\eqref{eq:ss_full_obj}.  The objective $f(x,w)$ is convex in $w$, but not strictly convex. 
To guarantee differentiability of $v(x)$, we add a smoothing term and consider 
\begin{equation}
\label{eq:vfbeta}
\begin{aligned}
v_\beta(x) := \min_{w} f(x,w) +  \frac{\beta}{2}\|w\|^2.
\end{aligned}
\end{equation}
where $\beta$ is usually taken to be a very small number (e.g. $10^{-4}$) so that the added term has minimal effect on the original value function. (The minimizer of $v_\beta$ is different from that of $v$.)
The function $v_\beta(x)$ is a \emph{Moreau envelope}~\cite[Def~1.22]{RockWets98} of the true value function $v$; we refer the interested reader to~\cite{aravkin2016variable} for details and examples concerning the Moreau envelope specifically (and nonsmooth variable projection more broadly).
The unique minimizer $w(x)$ can be found quickly and accurately since the minimization problem with respect to $w$ is strongly convex: 
projected gradient descent converges linearly and can be accelerated using the Fast Iterative Shrinkage-Thresholding Algorithm (FISTA)~\cite{beck2009fast} approach. With the minimizer $w(x)$, the gradient of $v_\beta$ is readily computed as
\begin{equation}
\label{eq:vgrad}
\nabla v_\beta(x) = \partial_x f(x,w)|_{w = w(x)}.
\end{equation}

Plugging $w(x)$ back into~\eqref{eq:ss_full_obj} we obtain the problem
\begin{equation}
\label{eq:ss_x_only}
\begin{aligned}
 \min_{x} v_\beta(x) = &\frac{1}{2} \sum_{t=0}^{T-1} \|\obst - \Meas(\cstatet)\|_{R^{-1}}^2 + \nu\|w_{t+1}(x) - w_{t}(x)\|_2^2 \\
 + & \sum_{m=1}^\ddim  w_{t,m}(x)r\log\left( 1 + \frac{\|\cstatetp -\Proc_m(\cstatet)\|_{Q^{-1}}^2}{r}\right) \\
 + &\nu\|w_{t+1}(x) - w_{t}(x)\|_2^2 + \frac{\beta}{2} \| w(x)\|^2,
\end{aligned}
\end{equation}
where $w_{t,m}(x) \equiv w_t [m](x)$. 


\subsection{Gauss-Newton descent with \Studentst penalties}
In this section we write the objective function \eqref{eq:ss_x_only} as a convex composite function and specify the line search we use to update $x$. The particular line search method was first proposed in \cite{burke1985descent} for a general class of algorithms for convex composite problems; 
the Gauss-Newton update derived here is a special case of this general line search scheme. 

To cast the objective into a convex composite function, let $v_\beta= \rho \circ F$, where
\[F(x) = \begin{pmatrix} f_1(x) \\ f_2(x) \end{pmatrix}\]
with
\begin{align*}
f_1(x) = &\frac{1}{2} \sum_{t=0}^{T-1} \sum_{m=1}^\ddim  w_{t,i}(x)r\log\left( 1 + \frac{\|\cstatetp -\Proc_m(\cstatet)\|_{Q^{-1}}^2}{r}\right) \\
&+ \nu\|w_{t+1}(x) - w_{t}(x)\|_2^2 + \frac{\beta}{2} \| w(x)\|^2 \\
 f_2(x) = & \Meas(x) - y
 \end{align*}
 and 
 \[\rho\begin{pmatrix} c \\ u \end{pmatrix} = c  + \frac{1}{2}\|u\|_{R^{-1}}^2 + \mnoise_{[0,+\infty]}(c).\]
 
At each iteration, we choose a search direction $d^*(x)$ that
\begin{equation}
\label{solve_d}
\begin{aligned}
d^* \in \text{argmin}_d & \quad  \rho(F(x)  +F^{(1)}(x)d) + \frac{1}{2}d^TU(x)d \\
\in \text{argmin}_d & \quad f_1(x) + \nabla f_1(x)d + \frac{1}{2}\|f_2(x) + \nabla f_2(x)d\|^2_{R^{-1}} \\
& \quad + \frac{1}{2}d^TU(x)d \\
\in \text{argmin}_d & \quad \frac{1}{2}d^T\left(U(x) + \nabla \Meas(x)^TR^{-1}\nabla \Meas(x)\right) d \\
& + \nabla v_\beta(x)^Td
\end{aligned}
\end{equation}
where the equivalence is obtained by dropping terms independent of $d$.
%
In general $U(x)$ can be any positive semidefinite matrix that varies continuously with respect to $x$, 
but for our particular objective function involving Student's t penalty, 
$U(x)$ is chosen to be a Hessian approximation of the \Studentst term in $f_1(x)$. Therefore the update can be interpreted as a Gauss-Newton style update. 
This approximation, proposed in~\cite[(5.5), (5.6)]{aravkin2014robust}, is employed here because of its significant computational advantage; it is
of the form
\begin{equation}
\label{u_matrix}
 U = \begin{bmatrix} 
U_1 & A_2^T & 0 &   \\
A_2 & U_2 & A_3^T & 0   \\
0 & \ddots & \ddots & \ddots \\
 & 0 & A_T & U_T 
\end{bmatrix} 
\end{equation}
with
\[A_t = - r \sum_{m=1}^\ddim w_{t-1,m}(x)\frac{Q^{-1}\nabla \Proc_m(x_{t-1})}{r + \|\cstatet - \Proc_m(x_{t-1})\|^2},\]
\begin{align*}
U_t = & r \sum_{m=1}^\ddim \frac{w_{t,m}(x)\nabla \Proc_m(\cstatet)^TQ^{-1}\nabla \Proc_m(\cstatet)}{r + \|\cstatetp-\Proc_m(\cstatet)\|^2} \\
+& \frac{w_{t-1,m}(x)Q^{-1}}{r + \|\cstatet - \Proc_m(x_{t-1})\|^2}
\end{align*}
for $1 \leq t \leq T-1$, and
\[ U_T = \frac{rw_{T-1,m}(x)Q^{-1}}{r + \|x_T - \Proc_m(x_{T-1})\|^2}.\]
We can rewrite $U(x)$ as
\[U(x)= \sum_m \Proc_m(x)^T \tilde{Q}_m(w(x))^{-1}\Proc_m(x),\]
where 
\[G_m(x) = 
\begin{bmatrix}
I & 0 & 0 &  \\
- \nabla \Proc_m(x_2) & I & 0 & 0 \\
0 & \ddots & \ddots & \ddots \\
\dots & 0 -\nabla \Proc_m(x_T) & & I 
 \end{bmatrix}
\]
and
\begin{align*}
 & \tilde{Q}_m(w(x))^{-1}  = \text{diag}(\tilde{Q}_{m,t}(w(x))^{-1}) \\
 & \tilde{Q}_{m,t}(w(x))^{-1} = \frac{r w_{t-1,m}(x)Q^{-1}}{r + \|x_{t}-\Proc_i(x_{t-1})\|^2}.
 \end{align*}
Clearly $U(x)$ is positive semidefinite; 
we show in Lemma~\ref{lemma:satisfied} that $U(x)$ is actually positive definite, 
so problem~\eqref{solve_d} reduces to the block triadiagonal linear system 
\[\left(U(x) +\nabla \Meas(x)^TR^{-1}\nabla \Meas(x)\right)d + \nabla v_\beta(x)  = 0. \]

Given $d^*(x)$, the new $x^+$ is of the form
\[x^+ = x + \mnoise d^*,\]
where $\mnoise$ is a step size selected using the \emph{Armijo}-type~\cite[Sec.~3.1]{nocedal2006nonlinear} line search criterion.
\begin{equation}
\begin{aligned}
\label{line_search}
\delta = \max \{\gamma^l\colon& \rho(F(x + \gamma^l d^*)) \leq \rho(F(x)) +  c\gamma^l \Delta(x;d^*) \\ 
&\text{and } c\in(0,1)\}
\end{aligned}
\end{equation}
%
with
\[\Delta(x;d)  = \rho(F(x) + F^{(1)}(x)d) + \frac{1}{2}d^TU(x)d - \rho(F(x)).\] 
When $d = 0$, we have $\Delta(x;0) = 0$%
\footnote{We overload $\Delta$ here to match the notation in \cite{burke1985descent,aravkin2014robust}; $\Delta(x;d^*)$ should not be confused with $\Delta^M$, which is used to denote the simplex containing relaxed state estimates.}%
, and since we choose the minimizing 
\[d^* = \argmin_d~\rho(F(x)  +F^{(1)}(x)d) + \frac{1}{2}d^TU(x)d,\]
we have 
$\Delta(x;d^*) \leq 0$.
Further, 
\begin{align*}
\Delta(x;d^*) = 0 &\Leftrightarrow 0 \in \argmin_d \rho(F(x)  +F^{(1)}(x)d) + \frac{1}{2}d^TU(x)d\\
& \Leftrightarrow 0 \in \partial \rho (F(x))F^{(1)}(x) 
\end{align*}
by \cite[Thm.~3.6]{burke1985descent}. In other words, stationarity is achieved when $\Delta(x;d^*) = 0$.
When $\Delta(x;d) < 0$, we are guaranteed to have descent
\[ \rho(F(x) + F^{(1)}(x)d) < \rho(F(x))\]
since $U(x)$ is positive semidefinite. 
%
This condition ensures that the line search step (\ref{line_search}) is well-defined \cite[Lemma~2.3]{burke1985descent}.

Our approach is summarized in Algorithm \ref{alg:ss_vp}. 
The positive parameter $\epsilon$
in the algorithm specifies the stopping condition.
Finally, 
we project the relaxed discrete state estimate ${w}_t \in \Delta^\ddim$ 
to obtain a discrete state estimate in $\mathcal{D}^\ddim$
as described in Section~\ref{sec:relax_smooth}.

\begin{algorithm}
  \caption[Caption]{\label{alg:ss_vp} Variable Projection for (\ref{eq:ss_full_obj}).} 
  \begin{algorithmic}[1]
    \Require{$x,w, Q,R,r,\nu,\beta,\epsilon$} 
    \For{$k = 1,2,3,... $} 
    	\Let{$d^{(k)}$}{Gauss-Newton direction for $x^{(k)}$} 
    	\label{alg:ss_vp:state_update_direction}
	\Let{$x^{(k+1)}$}{$x^{(k)} + \mnoise d^{(k)}$}
	\label{alg:ss_vp:state_update}
	\Let{$w^{(k+1)}$}{$\text{InnerSolver}_{\Pi_t\Delta}(w^{(k)})$}
	\Let{loss$_{k}$}{$f(x^{(k+1)},w^{(k+1)})$} \\
      \EndFor
      Iterate till $\Delta(x^{(k)};d^{(k)}) \geq -\epsilon$.
  \end{algorithmic}
  \label{alg:variable-projection}
\end{algorithm}

\subsection{Convergence of state estimation algorithm}

%
The convergence of Algorithm~\ref{alg:ss_vp} to a stationary point for a general class of convex composite objective functions is established in \cite{burke1985descent} and \cite{aravkin2014robust}. 
In particular ~\cite[Theorem~5.1]{aravkin2014robust} establishes the possible outcomes when applying this type of algorithm; informally, either the algorithm converges or the search direction $d_k$ diverges. 
In the remainder of this section we provide two technical results needed to formalize this intuition: 
%

%
 \begin{itemize}
\item Lemma~\ref{lemma:sufficient} establishes a set of sufficient conditions that prevent divergence ($\|d^{(k)}\|\to \infty$);
\item Lemma~\ref{lemma:satisfied} proves that the sufficient conditions are satisfied.
\end{itemize}


 \begin{lemma}
 \label{lemma:sufficient}
Let $\Lambda = \{y|\rho(y) \leq v_\beta(x^{(0)})\}$. If $F^{-1}(\Lambda) = \{x |F(x) \in \Lambda\}$ is bounded and $U(x)$ is positive definite for all $x \in F^{-1}(\Lambda)$,
then the hypotheses in \cite[Theorem~5.1]{aravkin2014robust} are satisfied and the sequence of search directions $\{d^{(k)}\}$ is bounded.
\end{lemma}

{\bf Proof}:
The hypotheses in \cite[Theorem~5.1]{aravkin2014robust} require that $F^{(1)}$ to be bounded and uniformly continuous on the set $S = \bar{co}(F^{(-1)}(\Lambda))$ where $\bar{co}$ stands for the closed convex hull. $F^{(1)}$ is continuous on $S$ since $f_1^{(1)}$ exists and is continuous by property of Moreau envelope and proximal operator, and $f_2^{(1)}$ is continuous trivially. Further, given that $S$ is closed by definition and bounded by assumption, it is compact. Hence $F^{(1)}$ is bounded and uniformly continuous on $S$.

Now we need to show that the sequence of search direction is bounded. At any iteration, the search direction $d$ we choose satisfies
\[ 0 \leq \rho(F(x) + F^{(1)}(x)d) + \frac{1}{2}d^TU(x)d \leq \rho(F(x)) \leq \rho(F(x^0))\]
where the first inequality relies on $\rho \geq 0$ and on the positive semidefinite property of $U(x)$; the second inequality comes from $\Delta(x;d) \leq 0$; the third inequality results from the line search condition that creates a decreasing sequence $\{\rho(F(x^{(k)})\}$.

Since $\rho(F(x^0))$ is finite, $d^TU(x)d < \infty$ for all iterations. Because $\Lambda$ is closed by closedness of $\rho$ and $F$ is continuous, $F^{-1}(\Lambda)$ is also closed. Along with its boundedness by assumption, $F^{-1}(\Lambda)$ is compact. Since $x \in F^{-1}(\Lambda)\mapsto\lambda_{\min}(U(x))$ is continuous, its image is bounded, hence given that $U(x)$ is positive definite there exists some $\lambda_{\min} > 0$ for all $x \in F^{-1}(\Lambda)$. Therefore $0 < \lambda_{\min}\|d\|^2 \leq d^TU(x)d < \infty$, which implies that $d^{(k)}$ cannot be unbounded. $\square$

\begin{lemma}
\label{lemma:satisfied}
$F^{-1}(\Lambda)$ is bounded for problem (\ref{eq:ss_x_only}) and $U(x)$ is positive definite for all $x \in F^{-1}(\Lambda)$.
\end{lemma}

{\bf Proof}:
First note that $\Lambda$ is bounded by the coercivity of $\rho$. This implies that for an unbounded sequence $\|x^{(k)}\|\to \infty$, we still have $f_1(x^{(k)}) < \infty$ and $\|f_2(x^{(k)})\| < \infty$.

If $\|x^{(k)}\|\to \infty$, then we can find some $t+1$ and a subsequence $J$ such that $\lim_{k \in J} \|\cstatetp^{(k)}\| = \infty$. By the definition of $f_1$ and $f_1(x^{(k)}) < \infty$, $\lim_{k \in J} \|\Proc_i(\cstatet^{(k)})\| = \infty$, which further implies that $\lim_{k \in J} \|x_{t}^{(k)}\| = \infty$. Iteratively this means that 
$\lim_{k \in J} \|x_{t}^{(k)}\| = \infty$ for all $t$, in particular for the given starting point $x_0$, but that is not possible.

To show that $U(x)$ in (\ref{u_matrix}) is positive definite, recall that we can rewrite $U(x)$ as 
\[U(x)= \sum_m G_m(x)^T \tilde{Q}_m(w(x))^{-1}G_m(x) \succeq 0.\]
If there exists some $d$ such that $d^TU(x)d = 0$, then 
\begin{align*}
 &d^T \left (\sum_m G_m(x)^T \tilde{Q}_m(w(x))^{-1}G_m(x)\right)d  \\
 = &\sum_m \underbrace{d^T G_m(x)^T}_{z_m(x)^T} \tilde{Q}_m(w(x))^{-1}\underbrace{G_m(x) d}_{z_m(x)}\\
 = & \sum_{m}z_m(x)^T\tilde{Q}_m(w(x))^{-1}z_m(x) = 0,\\
 \Rightarrow &z_m(x)^T\tilde{Q}_m(w(x))^{-1}z_m(x) = 0 ~\forall i \\
 \Rightarrow & z_{m,t}(x)^T \tilde{Q}_{m,t}(w(x))^{-1}z_{m,t}(x) = 0 \forall t~ \forall i
 \end{align*}
since $\tilde{Q}_m(w(x))^{-1}  = \text{diag}(\tilde{Q}_{m,t}(w(x))^{-1})$, and
\[\tilde{Q}_{m,t}(w(x))^{-1} = \frac{r w(x)_{t,m}Q^{-1}}{r + \|\cstatetp-\Proc_m(\cstatet)\|^2}\]
are positive semidefinite. However because each $w_t \in \Delta$, there has to be some $\tilde{Q}_{m,t}^{-1} \succ 0$ for each $t$. Therefore $U(x)$ must be positive definite for all $x \in F^{-1}(\Lambda)$. $\square$

\section{Experiments}
\begin{figure}
	\centering
	\begin{subfigure}[b]{.32\columnwidth}
	\centering
	\includegraphics[height=1.5in]{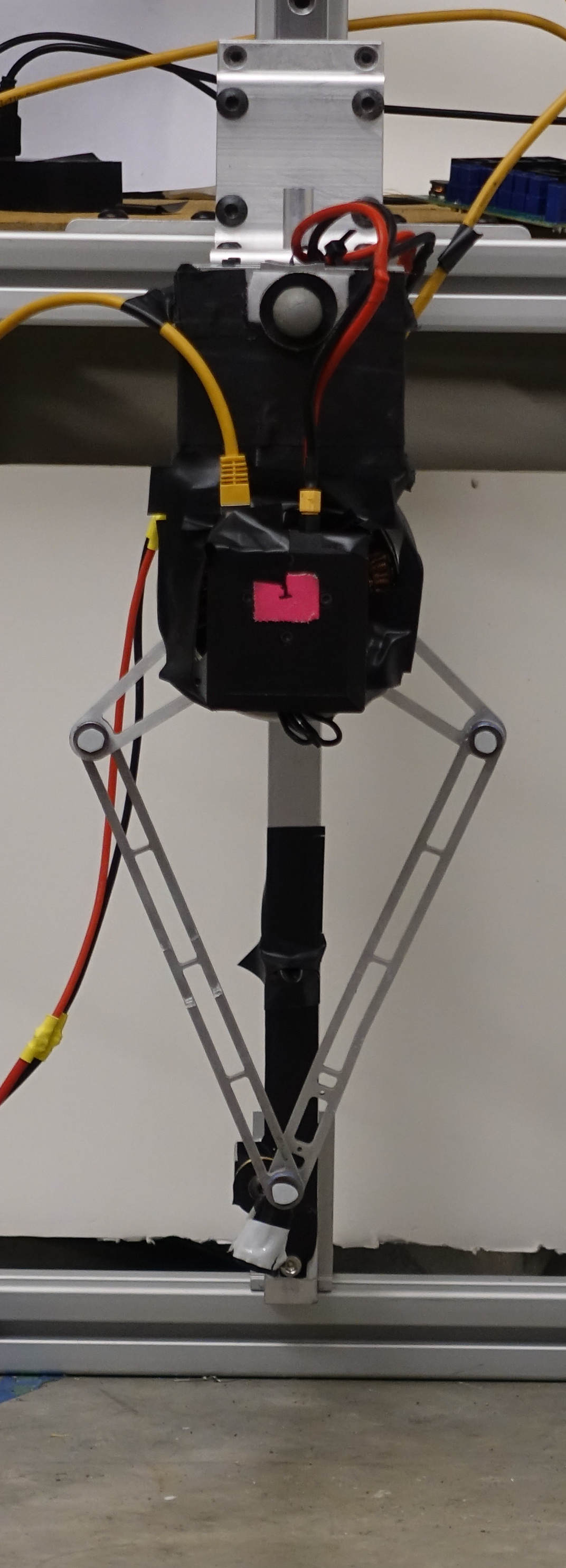}%
	\caption{}
	\label{fig:minitaur}
	\end{subfigure}
	\begin{subfigure}[b]{.32\columnwidth}
	\centering
	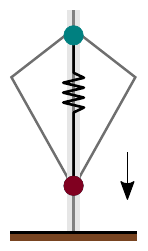%
	\caption{}
	\label{fig:nonlinear_hopper}
	\end{subfigure}%
	\begin{subfigure}[b]{.32\columnwidth}
	\centering
    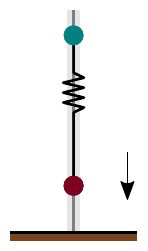%
	\caption{}
	\label{fig:linear_hopper}
	\end{subfigure}
	\caption{ \textbf{Jumping robot and impact oscillator hybrid system models} (Sec.~\ref{sec:dynamics_des}).
	(\subref{fig:minitaur}) Photograph of the physical robot (one leg from a Minitaur~\cite{KenneallyEtAl2016DesignPrinciplesFamily}) that inspired the simulation models.
    (\subref{fig:nonlinear_hopper}) Nonlinear model consisting of two masses coupled with a linear spring and a nonlinear pantograph mechanism.
    (\subref{fig:linear_hopper}) Linear model consisting of two masses coupled with a linear spring.
	}
	\label{fig:my_label}
\end{figure}

To evaluate the proposed approach to state estimation for hybrid systems, we apply our algorithm to linear and nonlinear impact oscillators. 
In addition to being well-studied~(\cite[\S1.2]{DiBernardo2008PiecewisesmoothDynamicalSystems},~\cite{Schatzman1998-gt}),
these mechanical systems were chosen since they are among the simplest physically-relevant models that satisfy our assumptions, including non--identity reset maps.
The parameter and trajectory regime considered in what follows is representative of a jumping robot constructed from one limb of a commercially-available quadrupedal robot~\cite{KenneallyEtAl2016DesignPrinciplesFamily} and controlled with an event-triggered stiffness adjustment;
\figref{fig:minitaur} contains a photograph of the limb.
The jumping robot's $\hip$ and foot are constrained to move vertically in a gravitational field, 
so the rigid pantograph mechanism 
depicted in~\figref{fig:nonlinear_hopper}
has two mechanical degrees-of-freedom (DOF)
coupled through nonlinear pin-joint constraints.
These two DOF are preserved, but their nonlinear coupling is neglected, in the piecewise-linear model illustrated in~\figref{fig:linear_hopper}.
The hybrid dynamics of these linear and nonlinear impact oscillators are specified in Section~\ref{sec:dynamics_des}
\label{sec:experiments}

We perform two sets of experiments.
The first set of experiments in Sec.~\ref{sec:exp_linear} 
concern the piecewise-linear model depicted in~\figref{fig:linear_hopper}
and
explore the consequences of our modelling assumptions and the efficacy of our proposed algorithm: 
\begin{itemize}
\item Sec.~\ref{sec:exp_students_t} demonstrates the advantage of employing a \Studentst distribution for process noise as compared to a Gaussian distribution;
\item Sec.~\ref{sec:exp_gauss_newton} demonstrates the superior convergence rate yielded by Gauss-Newton descent directions as compared to gradient (steepest) descent;
\item Sec.~\ref{sec:exp_hybrid_changes} demonstrates the advantage of 
smoothing the relaxed discrete state estimate;
and 
\item Sec.~\ref{sec:exp_linear_unobs_guard} demonstrates the algorithm's performance 
when \emph{onboard} measurements are used instead of \emph{offboard} measurements.
\end{itemize}
The second set of experiments in Sec.~\ref{sec:exp_nonlinear} evaluate our proposed approach using the nonlinear model depicted in~\figref{fig:nonlinear_hopper}.

Since this section is devoted to comparing estimated states to ground truth simulation results, and since our approach entails the determination of a relaxed discrete state estimate \emph{en route} to obtaining the discrete state estimate, we now introduce notation that distinguishes these quantities:
\begin{itemize}
    \item $\wtrue_t$ $\in \mathcal{D}^M$ denotes the ground truth discrete state;
    \item $\rwest_t$ $\in \Delta^M$ denotes the relaxed discrete state estimate;
    \item $\west_t$ $\in \mathcal{D}^M$ denotes the discrete state estimate.
\end{itemize}
This notational distinction was not introduced previously in the interest of readability 
since there was no ambiguity entailed by overloading notation in the problem formulation and algorithm specification. 

\subsection{Impact oscillator hybrid system models}
\label{sec:dynamics_des}
The continuous state 
$\cstate = (q, \dot{q})\in\R^4$
for the jumping robot hybrid system model
consists of the two-dimensional configuration vector $q\in\R^2$ and corresponding velocity $\dot{q}\in\R^2$,
where
$q[1]$ and $q[2]$ denote the vertical height of the $\hip$ and foot, respectively.
The foot is not permitted to penetrate the ground, 
$q[2] \ge 0$,
so the first part of the discrete state indicates whether this constraint is active:  $\A$ (air) if $q[2]>0$, $\G$ (ground) if $q[2]=0$.
To compensate for energy losses at impact, an event-triggered controller stiffens or softens a spring based on which direction the hip is traveling, so the second part of the discrete state indicates the direction of travel for $q[1]$:
$\uparrow$ if up,
$\downarrow$ if down.
With 
$\ddot{q}_m(q,\dot{q})\in\R^2$ 
denoting the acceleration of the hip and foot in discrete state $m\in\set{\AD,\GD,\GR,\AR}$, 
formulae for this acceleration are given in Table~\ref{table:hopper_modes}.
At the moment of impact (when the discrete state changes from $\delem_t\in\{\AD,\AR\}$ to $\delem_{t+1}\in\{\GD,\GR\}$)
the foot velocity $\dot{q}[2]$ is instantaneously \emph{reset} to $0$, corresponding to perfectly \emph{plastic} impact.
An example of the jump in continuous state when transitioning from $\AD$ to $\GD$ on the foot velocity $\dot{q}[2]$ is shown in \figref{fig:ss_st} near time 17.5s.

\newcolumntype{C}{ >{\centering\arraybackslash} m{1.1 cm} }
\begin{table}[h]
\centering
\begin{tabular}{c C c}
\textbf{Discrete state $w$} & \textbf{Icon} & $\ddot{q}_\delem(x)$ \\
\hline
$\delem=\AD$ &  \includegraphics{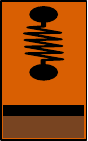} &  
	$\begin{bmatrix}
	\frac{1}{\mh} \left( -k_1(q,\dot{q})\right) -g  \\
	\frac{1}{m_t} \left( k_1(q, \dot{q})\right) -g  \\
	\end{bmatrix} $  \\
$\delem=\GD$ &  \includegraphics{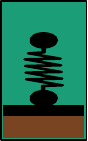} &  
	$\begin{bmatrix}
	\frac{1}{\mh} \left( -k_1(q,\dot{q})\right) -g  \\
	0 \\
	\end{bmatrix} $  \\
$\delem=\GR$ &  \includegraphics{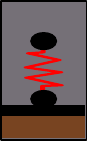} &  
	$\begin{bmatrix}
	\frac{1}{\mh} \left( -k_2(q,\dot{q})\right) -g  \\
	0 \\
	\end{bmatrix} $  \\
$\delem=\AR$ &  \includegraphics{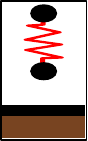}  &
	$\begin{bmatrix}
	\frac{1}{\mh} \left( -k_2(q,\dot{q})\right) -g  \\
	\frac{1}{m_t} \left( k_2(q, \dot{q})\right) -g  \\
	\end{bmatrix} $  \\
\end{tabular}
\caption{\textbf{Discrete states and continuous dynamics for impact oscillator hybrid system models} (Sec.~\ref{sec:dynamics_des}).  
Note that the continuous dynamics $\ddot{q}$ have the same general form for both the piecewise-linear and -nonlinear models, with the spring law $k$ being a linear or nonlinear function of the continuous state $x = (q,\dot{q})$ depending on which model is considered.
}
\label{table:hopper_modes}
\end{table}

\subsection{Piecewise-linear impact oscillator experiment}
\label{sec:exp_linear}
In this subsection, we employ the linear spring laws 
\[k_1(q,\dot{q}) = 10(q[1]-q[2])-3, \]
\[k_2(q,\dot{q}) = 15(q[1]-q[2])-3,\]
with parameter values $\mh = 3, m_t = 1, g=2$.

%

In our first demonstration 
the observed states are $q[1]$ and $q[2]$, position of the $\hip$ and foot, leaving the velocities unobserved:
\begin{equation}
\Meas_{\text{pos}}(x) = q.
\label{eq:H_pos}
\end{equation}
State estimation results for this system are shown in \figref{ss_xw_sm}.

In the remainder of this subsection, we demonstrate the effects of 
the choices we made in our problem formulation (Sec.~\ref{sec:problem}) and algorithm derivation (Sec.~\ref{sec:algorithm}) using the piecewise-linear model as a running example. 
We also consider a variation where the measurements correspond to the leg length and velocity,
which are more representative of the \emph{onboard} measurements available to an autonomous robot operating outside of the laboratory.
%


\subsubsection{\Studentst versus Gaussian process noise}%
\begin{figure}
\centering
\includegraphics[scale = 1]{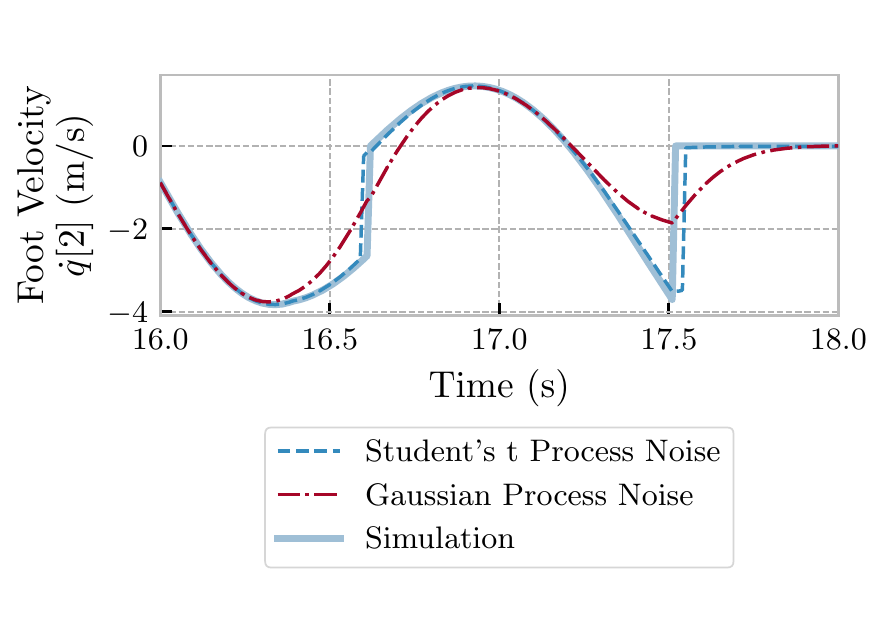}
\caption{\textbf{The \Studentst distribution process noise yields better estimates of instantaneous changes in continuous \state} (Sec.~\ref{sec:exp_students_t}).
In this plot, estimates of the foot velocity are shown
near two impacts ($\approx 16.6$s, $17.5$s). 
}%
\label{fig:ss_st}
\end{figure}

\label{sec:exp_students_t}
\figref{fig:ss_st} compares the estimation of foot velocity using \Studentst with $r = 0.01$ versus using Gaussian for the process noise distribution; in both cases the true discrete state is given.
The estimated trajectory for both distributions match the true simulated trajectory
away from jumps, while near jumps, such as around times $16.6$s and $17.5$s, using the \Studentst distribution enables closer tracking of the instantaneous change in the true foot velocity $\dot{q}[2]$ than when using a Gaussian distribution.

\subsubsection{Gauss-Newton versus gradient (steepest) descent }%
\begin{figure}
\centering
\includegraphics[scale = 1]{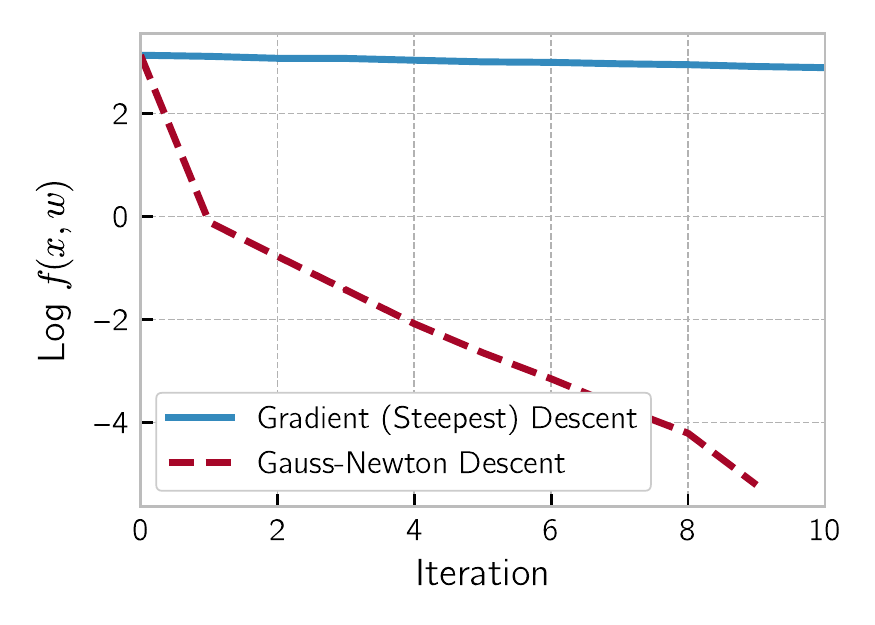}
\caption{
\label{fig:gn_gd}
\textbf{Gauss-Newton descent directions yield faster
convergence than gradient (steepest) descent} (Sec.~\ref{sec:exp_gauss_newton}).
In this plot, the discrete state variables $w$ are given and the second line of Algorithm~\ref{alg:variable-projection} is modified to use either Gauss-Newton descent directions or gradient (steepest) descent to estimate the continuous state variables $x$ by minimizing the relaxed objective function $f(x,w)$~\eqref{eq:ss_full_obj}.}
\end{figure}

\label{sec:exp_gauss_newton}
We empirically compared convergence rates for continuous state $\cstatet$ updates obtained using Gauss-Newton and gradient (steepest) descent directions 
(Algorithm \ref{alg:variable-projection}, line \ref{alg:ss_vp:state_update_direction}).
\figref{fig:gn_gd} shows the log loss versus algorithm iteration for the two methods; the actual discrete state $w_t$ was taken as given to perform this comparison. 
As expected, the objective value decreases significantly faster
when the search direction is determined by the Gauss-Newton scheme as compared to the direction of steepest descent,
reaching the stopping criterion in ten times fewer iterations in our tests.

\subsubsection{Smoothing the relaxed discrete state versus not}
\label{sec:exp_hybrid_changes}
\begin{figure}
\begin{center}  
\scalebox{1.0}{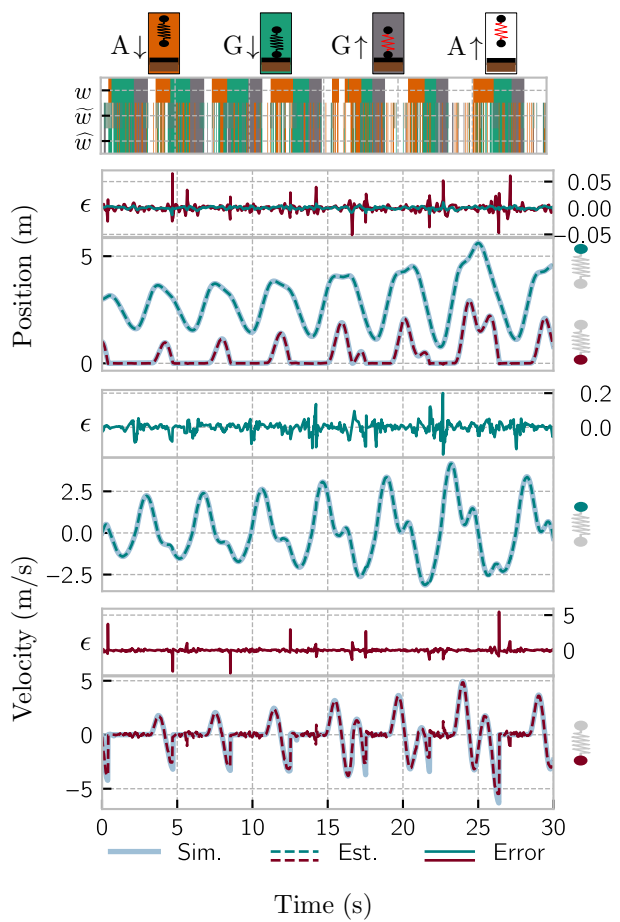}
\end{center}
\caption{\textbf{Without smoothing ($\nu = 0$), the discrete state estimate switches frequently} (Sec.~\ref{sec:exp_hybrid_changes}).
The top plot shows 
the true discrete state of the system $w\in D^M$, 
the relaxed discrete state estimate $\widetilde{w}\in \Delta^M$,
and 
the discrete state estimate $\widehat{w}\in D^M$
for a simulation of the piecewise-linear system.
The subsequent plots show the estimate, simulation, and error $\epsilon$ values for position and velocity of the hip and foot.}
\label{ss_xw_nsm}
\end{figure}

\begin{figure}
\begin{center}  
\scalebox{1.0}{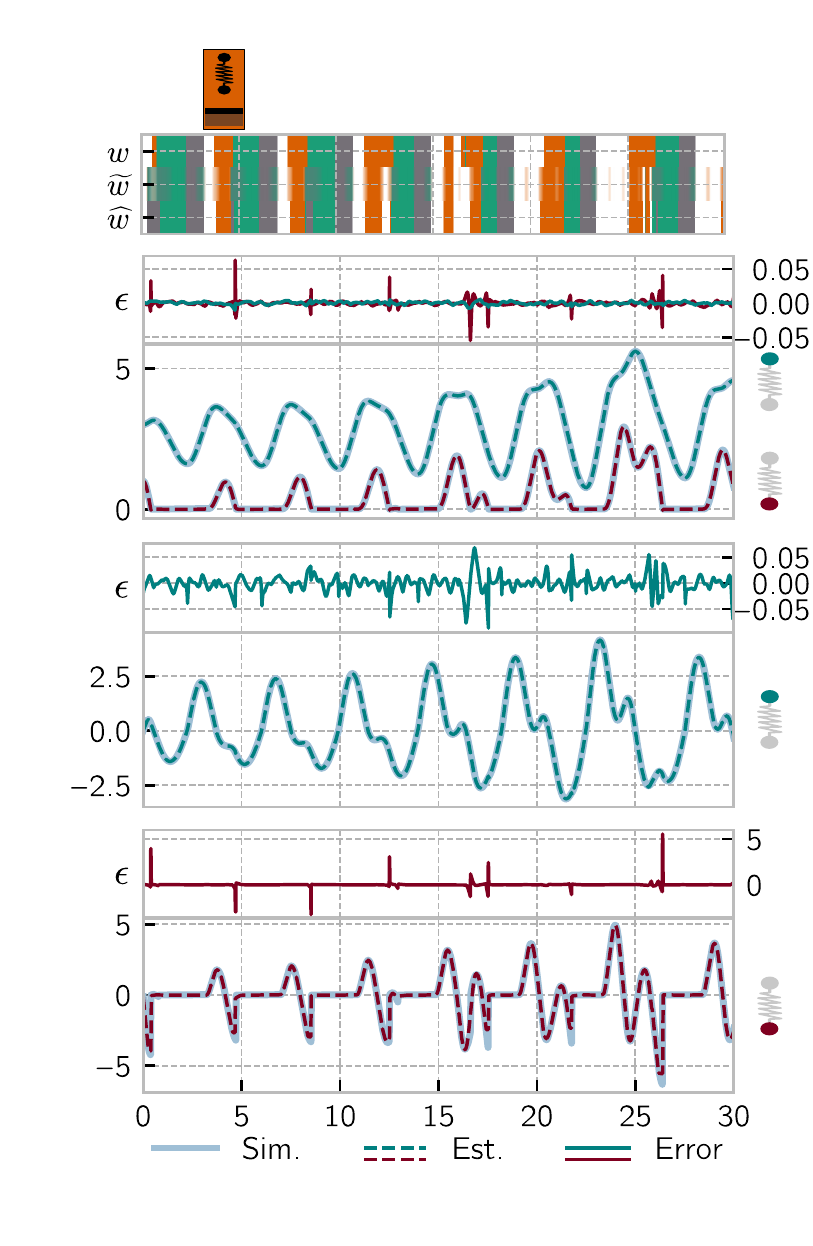}
\end{center}
\caption[With smoothing ($\nu > 0$), the estimated discrete state closely follows true discrete state.]{\textbf{With smoothing ($\nu > 0$), the discrete state estimate mostly switches near the true switching times.} (Sec.~\ref{sec:exp_hybrid_changes}).
This plot shows results from the piecewise-linear system;
the notational and plotting conventions are adopted from~\figref{ss_xw_nsm}.
}
\label{ss_xw_sm}
\end{figure}

If the continuous states are given, 
the discrete state estimate returned by our algorithm 
(skipping lines 2 and 3 of Algorithm~\ref{alg:variable-projection}) 
is very close to the true discrete state regardless of whether a smoothing term is included in the relaxed problem formulation.
When simultaneously estimating both the continuous and discrete states, the smoothing term becomes crucial,
as illustrated by comparing the discrete state estimates ($\widehat{w}_t$) in
\figref{ss_xw_nsm} (without smoothing) 
and \figref{ss_xw_sm} (with smoothing). 
In particular, 
the estimated discrete state switches rapidly without smoothing, 
whereas with smoothing the discrete state tends to remain constant for many samples and change mostly near ground-truth switching times.

\subsubsection{Onboard versus offboard measurements}
\label{sec:exp_linear_unobs_guard}
\begin{figure}
\centering
\scalebox{1.0}{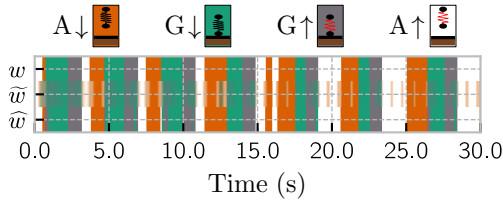}
\caption{\textbf{Estimated discrete state 
using onboard (relative position and velocity) measurements $\Meas_{\text{relative}}$~\eqref{eq:H_relative} for the piecewise-linear system closely matches true discrete state. 
} (Sec.~\ref{sec:exp_linear_unobs_guard}).
Continuous state estimates are not shown since they are formally unobservable using only onboard measurements (in practice, they drift away from ground truth over time).}
\label{ss_angle_nr}
\end{figure}

In the laboratory, the positions of the robot hip and foot can be directly measured \emph{offboard}, e.g. with an external camera system.
Outside of the laboratory, only the relative position of the hip and foot can be directly measured onboard our robot. 
Thus, we are motivated by this practical consideration to evaluate our algorithm's performance in the case where only the relative position and velocity of the hip and foot are measured,
\begin{equation}
\Meas_{\text{relative}}(x) = \begin{bmatrix} q[1] - q[2] \\ \dot{q}[1] - \dot{q}[2]\end{bmatrix}.
\label{eq:H_relative}
\end{equation}
Although the full hybrid system state is formally unobservable with these relative measurements,
our algorithm nevertheless yields good estimates of the discrete state as shown in~\figref{ss_angle_nr}; due to large errors in the estimate of (unobservable) continuous states, we omit those results from the figure.
%

\subsection{Piecewise-nonlinear impact oscillator experiment}
\label{sec:exp_nonlinear}
\begin{figure}
\centering
\scalebox{1.0}{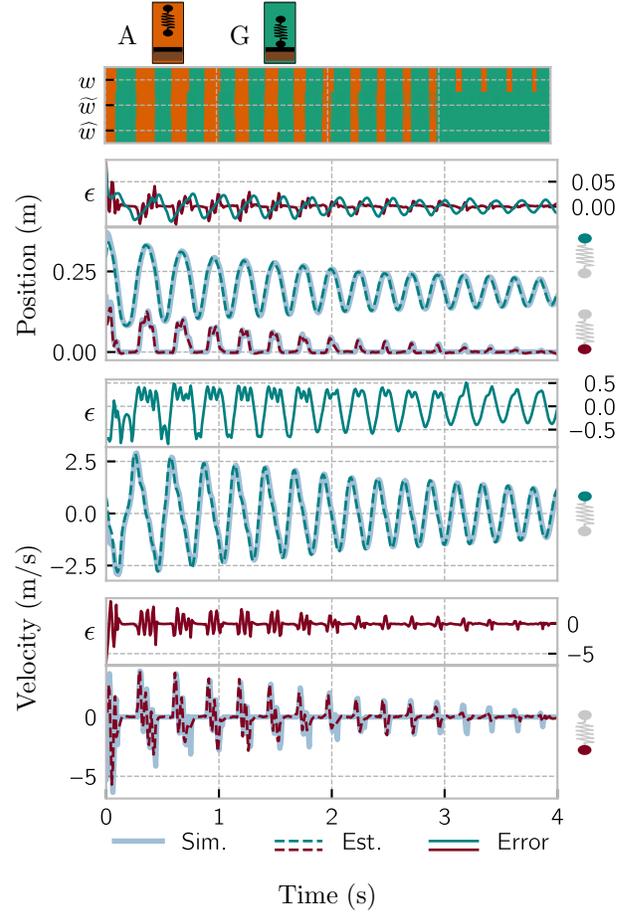}
\caption{\textbf{Continuous and discrete states estimated for the piecewise-nonlinear model} (Sec.~\ref{sec:exp_nonlinear}).
Notational and plotting conventions are adopted from~\figref{ss_xw_nsm};
note that this model only has two discrete states (Sec.~\ref{sec:dynamics_des}).
}
\label{fig:nonlinear}
\end{figure}

\begin{figure}
\centering
\scalebox{1.0}{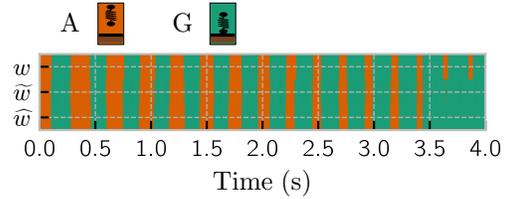}
\caption{\textbf{Estimated discrete state 
using onboard (relative position and velocity) measurements $\Meas_{\text{relative}}$~\eqref{eq:H_relative} for the piecewise-nonlinear system closely matches true discrete state.
} (Sec.~\ref{sec:exp_linear_unobs_guard}).
As with~\figref{ss_angle_nr}, continuous state estimates are not shown since they drift from the true values over time;
note that this nonlinear model only has two discrete states (Sec.~\ref{sec:dynamics_des}).
}
\label{fig:nonlinear:relative}
\end{figure}

To test Algorithm \ref{alg:variable-projection} on a nonlinear model, 
we included the kinematic constraints depicted in \figref{fig:nonlinear_hopper}, resulting
in a nonlinear spring force. 
In this model we set the two spring laws to be the same $k_1 = k_2$, decreasing the number of discrete states from four to two: 
$\delem = \A$ when $q[2]>0$ and $\delem = \G$ when $q[2]=0$.
State estimation results 
compare favorably with the analogous results from the piecewise-linear system when using either absolute position measurements $\Meas_{\text{pos}}$~\eqref{eq:H_pos} (compare~\figref{fig:nonlinear} with \figref{ss_xw_sm})
or
relative measurements $\Meas_{\text{relative}}$~\eqref{eq:H_relative}
(compare~\figref{fig:nonlinear:relative} with~\figref{ss_angle_nr}).

In \figref{fig:nonlinear} we see that the model can estimate continuous and discrete states in the nonlinear setting. 
However, we do notice that the estimated trajectories are not as close to ground truth as in the linear case. 
In particular, when 
$q[2]$
has a value only slightly greater than 0 (e.g. between times 3s and 4s), 
the algorithm fails to detect the transition between $\delem=\A$ and $\delem=\G$.

\section{Conclusion}
We proposed a new state estimation algorithm for hybrid systems, analyzed its convergence properties, and evaluated its performance on piecewise-linear and -nonlinear hybrid systems with non-identity resets.
The algorithm leverages a relaxed  state estimation problem formulation where the decision variables corresponding to the discrete state are allowed to take on continuous values.
This relaxation yields a continuous optimization problem that can be solved using recently-developed nonsmooth variable projection techniques.
The effectiveness of the approach was demonstrated on hybrid system models of mechanical systems undergoing impact.

\appendix

\begin{ack}                               
\end{ack}

\bibliographystyle{plain}        
\bibliography{refs}           
\end{document}